\documentclass[preprint,12pt]{elsarticle}
\usepackage{graphicx} 
\usepackage{amssymb}
\usepackage{amsthm}
\usepackage{graphics}
\usepackage{epsfig}
\usepackage{amssymb}
\usepackage{hyperref}
\usepackage{graphicx}
\usepackage{subfigure}
\usepackage{float}
\usepackage{color}
\usepackage{tikz-network}
\usepackage{tikz}
\usepackage{tikz,pgfplots}
\usetikzlibrary{patterns}
\usepackage[ruled,vlined]{algorithm2e}
\usepackage{hyperref}
\usepackage{caption}
\DeclareCaptionLabelFormat{cont}{#1~#2\alph{ContinuedFloat}}
\usepackage{a4wide}
\usepackage{amsmath}
\usepackage{amsfonts}
\usepackage{enumerate}
\usepackage{tikz}
\newcommand{\pf}{\noindent {\bf Proof: }}

\newtheorem*{theoremaux}{Theorem \theoremauxnum}
\gdef\theoremauxnum{1}

\newtheorem{lemma}{\bf Lemma}[section]
 
\newtheorem{theorem}{\bf Theorem}[section]
\newtheorem{conjecture}{\bf Conjecture}[section]
\newtheorem{proposition}[lemma]{\bf Proposition}
\newtheorem{corollary}[lemma]{\bf Corollary}
\newtheorem{definition}{\bf Definition}[section]
\newtheorem{remark}{\bf Remark}[section]

\journal{~}

\begin{document}

\begin{frontmatter}
\title{{Non-isomorphic $d$-integral circulant graphs}}





\author[1]{Sauvik Poddar}
\ead{sauvikpoddar1997@gmail.com}

\author[1]{Angsuman Das\corref{cor1}}
\ead{angsuman.maths@presiuniv.ac.in}

\address[1]{Department of Mathematics, Presidency University, Kolkata, India\\
86/1, College Street, Kolkata, India.
}
\cortext[cor1]{Corresponding author}

\begin{abstract}
The algebraic degree $Deg(G)$ of a graph $G$ is the dimension of the splitting field of the adjacency polynomial of $G$ over the field $\mathbb{Q}$. It can be shown that for every positive integer $d$, there exists a circulant graph with algebraic degree $d$. Let $C(d)$ be the least positive integer such that there exists a circulant graph of order $C(d)$ having algebraic degree $d$. A graph $G$ is called $d$-integral if $Deg(G)=d$. We call a $d$-integral circulant graph \textit{minimal} if order of that graph equals $C(d)$. Let $\mathcal{F}_{n,d}$ denote the collection of isomorphism classes of connected, $d$-integral circulant graphs of some given possible order $n$. In this paper we compute the exact value of $C(d)$ and provide some bounds on $|\mathcal{F}_{n,d}|$, thereby showing that the minimal $d$-integral circulant graph is not unique. Moreover, we find the exact value of $|\mathcal{F}_{p,d}|$ where both $p$ and $d$ are prime.
\end{abstract}

\begin{keyword}
	algebraic degree \sep Cayley graph \sep integral graph
	\MSC[2008] 05E18, 05C25, 05C30
	
\end{keyword}
\end{frontmatter}


\section{Introduction}

The notion of algebraic degree is introduced by M{\"o}nius et al., 2018 \cite{monius2018graphs},  as a generalization of the integral graphs i.e., graphs whose spectrum consists entirely of integers. The quest of characterizing integral graphs was first initiated by Harary and Schwenk \cite{harary1974graphs} in their $1973$ paper. However this project was later realized to be quite intractable as the author themselves mentioned so in their paper. Nevertheless this motivated the authors in \cite{monius2018graphs} to pose the following notion of algebraic degree of a graph. 

\begin{definition}\cite{monius2018graphs}
Given a graph $G$, the algebraic degree, $Deg(G)$, is the dimension of the
splitting field of the adjacency polynomial of $G$ over the field $\mathbb{Q}$. 
\end{definition}

\begin{definition}\cite{abdollahi20242}
A graph $G$ is called $d$-integral if $Deg(G)=d$. 
\end{definition}

Since the eigenvalues of a graph are algebraic integers, a graph $G$ is integral if and only if $Deg(G)=1$. The splitting field of the adjacency polynomial (i.e., the characteristic polynomial of the adjacency matrix) of a graph $G$ is called the \textit{splitting field} of $G$ \cite{monius2022splitting}. We denote the splitting field of a graph $G$ by $\mathbb{SF}(G)$. Hence, $Deg(G)=[\mathbb{SF}(G):\mathbb{Q}]$.

In \cite{monius2018graphs} the authors investigated the graph properties for which the spectrum becomes non-integral. They found an interesting result that a graph with sufficiently large diameter has some eigenvalues of large degree. In other way around, it can be interpreted as given a graph with some fixed algebraic degree and fixed maximum vertex degree, its diameter must be bounded above. Loosely speaking, algebraic degree of a graph is a parameter which measures how far a graph is from being an integral one.


Let $G=(V,E)$ be a simple, undirected graph. By $n_G$, we denote the order of $G$ and by $sp(G)$, we denote the spectrum (i.e., the multiset of the adjacency eigenvalues) of $G$. Two graphs $G$ and $H$ are said to be \textit{isospectral} (or \textit{cospectral}) if $sp(G)=sp(H)$. For graphs $G$ and $H$, $G\cong H$ implies that $sp(G)=sp(H)$. In this paper, we deal with circulant $d$-integral graphs. 
\begin{definition}
A circulant graph $\Gamma=Cay(\mathbb{Z}_n,S)$ is a Cayley graph over $\mathbb{Z}_n$ with respect to an inverse-symmetric subset $S\subseteq\mathbb{Z}_n$. The set $S$ is also called the connection set of $\Gamma$. 
\end{definition}
The cycle graph $\mathcal{C}_n\cong Cay(\mathbb{Z}_n,\lbrace{1,-1}\rbrace)$ and the complete graph $\mathcal{K}_n\cong Cay(\mathbb{Z}_n,\mathbb{Z}_n\setminus \{0\})$ are some well-known examples of circulant graphs. Since we consider loopless, undirected circulant graphs, we take $0\notin S$ and $S=-S$. 

Previous works on algebraic degree of circulant graphs can be found in \cite{monius2020algebraic} and \cite{monius2022splitting}. For a subgroup $H\le Aut(\mathbb{Z}_n)$, let $\mathbb{Q}(\zeta_n)^{H}$ denote the fixed field of $H$, i.e. the unique maximum
subfield of $\mathbb{Q}(\zeta_n)$ where each element is fixed by every automorphism in H. Moreover, for an element $x\in\mathbb{Z}_n$, let $x^{H}:=\lbrace{\sigma(x):\sigma\in H}\rbrace$ be the orbit of $x$ under $H$, and for $S\subseteq \mathbb{Z}_n$ let $S^{H}:=\bigcup_{s\in S}s^{H}$. Then $S$ is said to be fixed by $H$ if $S^{H}=S$. (These notations and terminologies are adapted from \cite{monius2022splitting}). If we identify $Aut(\mathbb{Z}_n)\cong\mathbb{Z}_{n}^{*}$, then the orbit of $x$ under $H$ is nothing but the left coset $xH$ of $H$ in $\mathbb{Z}_{n}^{*}$, i.e., for $S\subseteq \mathbb{Z}_n$ the set $S^{H}$ becomes $S^{H}:=\bigcup_{s\in S}sH$. The main results of \cite{monius2020algebraic} and \cite{monius2022splitting} are Theorem \ref{thm-circulant-alg-deg-d} and Theorem \ref{split-field-circulant}, which will be used later.

Another area of research involves the classification of the isomorphism classes of circulant graphs of order $n$. Although the complete classification of the isomorphism classes of circulant graphs of arbitrary order has not been resolved yet, several works on this topic has been done extensively over the years. For example, the complete classification of the isomorphism classes of integral circulant graphs of order $n\in\mathbb{N}$ has been resolved by Klin and Kov{\'a}cs \cite{klin2012automorphism}. The complete classification of the isomorphism classes of circulant graphs of order $p^2$ $(p>2)$ was done by Alspach and Parsons \cite{alspach1979isomorphism}, of order $p^k$ $(p>2,k\ge1)$ was done by Klin and P{\"o}schel \cite{klin1980isomorphism}, and of order $2^k$ $(k\ge1)$ was done by Klin, Najmark and P{\"o}schel \cite{klin1981Schur}.

\subsection{\textbf{Our Contribution}}
In this present article, we address two closely related issues in connection with algebraic degree of circulant graphs. Firstly, in Section \ref{ord-min-d-int-circ-graph}, for a given positive integer $d$, we compute the minimum order of a $d$-integral circulant graph by providing an explicit construction for the same. To be more precise, we want to find the least positive integer $t$ such that there exists a circulant graph on $t$ vertices having algebraic degree $d$. We denote this $t$ (which is indeed a function of $d$) by $C(d)$. We refer such graphs as \textit{minimal} $d$-integral circulant graphs. To study $C(d)$, we first ensure the existence of a $d$-integral circulant graph for any given positive integer $d$ in Theorem \ref{theorem-circulant-alg-deg-d}. In Theorem \ref{C(d)}, we provide the exact value of $C(d)$ and compare it with the bound obtained from Corollary \ref{cor-circ-min-ord} in Table \ref{tab:my_label}. 

Secondly, in Section \ref{count-number-conn-d-int-circ-graphs}, in view of the works done in \cite{alspach1979isomorphism},\cite{klin2012automorphism},\cite{klin1981Schur} and \cite{klin1980isomorphism}, we focus on finding the number of isomorphism classes of connected, $d$-integral circulant graphs of some given possible order. More precisely, we determine the exact number of connected, non-isomorphic, integral circulant graphs of order $n\in\mathbb{N}$. Also for a given integer $d>1$, we provide some lower and upper bounds on the number of connected, non-isomorphic, $d$-integral circulant graphs of some given possible order. As a result, for a prime $p>3$, we establish a lower bound regarding the number of connected, non-isomorphic circulant graphs on $p$ vertices with maximum algebraic degree, which is a betterment of Corollary 2.8 in \cite{monius2020algebraic}. Finally for a prime $d$, we determine the exact number of non-isomorphic, $d$-integral circulant graphs of a given prime order. 

\section{Order of a minimal $d$-integral circulant graph}\label{ord-min-d-int-circ-graph}
In this section, we show that for every positive integer $d$, there exists a graph with algebraic degree $d$. Though the existence was implicitly mentioned in \cite{monius2022splitting}, we provide a formal and constructive proof of the same. To begin with, we recall a result from \cite{monius2020algebraic}.

\begin{theorem}\label{thm-circulant-alg-deg-d}[\cite{monius2020algebraic}, Theorem $2.5$]
Let $p>2$ be a prime and $S\subseteq\mathbb{Z}_p$ be a connection set of the circulant graph $\Gamma:=Cay(\mathbb{Z}_p,S)$ with $|S|=M$. Let $m$ be the maximum common divisor of $M$ and $p-1$ such that $S$ is $m$-decomposable i.e., $S$ can be written as $$S=\bigcup_{i=1}^{M/m}S_i$$ with $|S_i|=m$ such that $s_{i,1}^m\equiv\cdots\equiv s_{i,m}^m~(mod~p)$ for $s_{i,j}\in S_i$, $i=1,\ldots,M/m$ and $j=1,\ldots,m$.
Then, $$Deg(\Gamma)=\frac{p-1}{m}.$$ 
\end{theorem}

\begin{theorem}\label{theorem-circulant-alg-deg-d}
For any positive integer $d$, there exists a circulant graph of prime order with algebraic degree $d$.
\end{theorem}
\pf Let $d\in\mathbb{N}$ be fixed. We choose a prime $p$ of the form $2dk+1$, where $k\in\mathbb{N}$. Set $m=2k$ so that $p=dm+1$. Let $r$ be a primitive root in $\mathbb{Z}_p^{*}$. Consider the set $S=\lbrace{1,r^d,r^{2d},\ldots,r^{(m-1)d}}\rbrace$. Then $|S|=m=M$ and $gcd(M,p-1)=m$. It is easy to verify that $S$ is an inverse-symmetric subset of $\mathbb{Z}_p$. Since $$(r^{jd})^m\equiv(r^{md})^j\equiv(r^{p-1})^j\equiv1~(mod~p)$$ for all $0\le j\le m-1$, it follows that $S$ is $m$-decomposable. Hence by Theorem \ref{thm-circulant-alg-deg-d}, $Deg(Cay(\mathbb{Z}_p,S))=\frac{p-1}{m}=d$.\qed

In view of the above theorem, we define $$C(d)=\min \{n_{\Gamma}:\Gamma \mbox{ is circulant and } Deg(\Gamma)=d\}.$$ Thus we have the following corollary.

\begin{corollary}\label{cor-circ-min-ord}
For $d\in\mathbb{N}$, let $p_d$ denote the smallest prime of the form $p_d\equiv 1~(mod~2d)$. Then $C(d)\le p_d$.
\end{corollary}

Now we focus on the order of a minimal $d$-integral circulant graph. In \cite{monius2022splitting}, M{\"o}nius proved the following result.

\begin{theorem}\label{split-field-circulant}[\cite{monius2022splitting},Theorem 5.3]
The splitting field of $Cay(\mathbb{Z}_n,S)$ is given by $\mathbb{Q}(\zeta_n)^{H}$, where $H$ is the maximum subgroup of $Aut(\mathbb{Z}_n)\cong\mathbb{Z}_n^{*}$ such that $S^{H}=S$.
\end{theorem}
Let $\varphi(\cdot)$ denote the Euler's totient function. Then the algebraic degree of a circulant graph of order $n$ is given by the following expression, which is Theorem $6.1$ in \cite{monius2022splitting}.

\begin{equation}\label{formula-alg-deg-circ}
Deg(Cay(\mathbb{Z}_n,S))=\dfrac{\varphi(n)}{|H|}
\end{equation}
where $H=\lbrace{k\in\mathbb{Z}_n^{*}:kS=S}\rbrace$ is the maximum subgroup of $\mathbb{Z}_n^{*}$ fixing $S$.

It is easy to observe that any subgroup of $\mathbb{Z}_n^{*}$ which fixes 
$S$ must contain $1$ and $-1$. Therefore any such $H$ must be of even order and thus the algebraic degree of the family of circulant graphs $Cay(\mathbb{Z}_n,S)$ divides $\varphi(n)/2$. Interestingly the converse is also true, which we state below (Proposition \ref{divisor-phi-circulant-prop}), preceded by a lemma which is used in proving the result.

\begin{lemma}\label{divisor-phi-circulant-lemma}
Let $n\ge 3$ and $l$ be an even divisor of $\varphi(n)$. Then there exists a subgroup of $\mathbb{Z}_n^{*}$ of order $l$ which is inverse-symmetric in $\mathbb{Z}_n$.
\end{lemma}
\pf For an even divisor $l$ of $\varphi(n)$, choose a subgroup $H$ of $\mathbb{Z}_n^{*}$ of order $l$. If $H$ is inverse-symmetric in the group $\mathbb{Z}_n$, then we are done. If not, then take $K\le H$ such that $[H:K]=2$. Let $T:=K\cup -K$. It is easy to check that $K\cap -K=\emptyset$ and $T$ is a subgroup of $\mathbb{Z}_n^{*}$. Also $T$ is inverse-symmetric in $\mathbb{Z}_n$ and $|T|=2|K|=l$.\qed

\begin{proposition}\label{divisor-phi-circulant-prop}
Let $n\ge 3$. Then for every divisor $d$ of $\frac{\varphi(n)}{2}$, there is a $\frac{\varphi(n)}{d}$-regular circulant graph on $n$ vertices 
of algebraic degree $d$.
\end{proposition}

\pf Let $d\mid \frac{\varphi(n)}{2}$. Then $\frac{\varphi(n)}{d}$ is an even divisor of $\varphi(n)$. By Lemma \ref{divisor-phi-circulant-lemma}, there exists a subgroup, say $S$ of $\mathbb{Z}_n^{*}$ such that $|S|=\frac{\varphi(n)}{d}$, which is inverse-symmetric in $\mathbb{Z}_n$. Let $\Gamma=Cay(\mathbb{Z}_n,S)$. Since $S$ is the largest subgroup of $\mathbb{Z}_n^{*}$ fixing $S$ itself, by Theorem \ref{split-field-circulant}, $\mathbb{SF}(\Gamma)=\mathbb{Q}(\zeta_n)^{S}$ and hence $Deg(\Gamma)=\frac{\varphi(n)}{|S|}=d$ by \ref{formula-alg-deg-circ}. Also $\Gamma$ is $|S|=\frac{\varphi(n)}{d}$ regular. This completes the proof.\qed


For $d\in\mathbb{N}$, consider the set $\mathcal{N}_d:=\lbrace{n\in\mathbb{N}:2d\mid\varphi(n)}\rbrace$. It is easy to see that for all $d\in\mathbb{N}$, the set $\mathcal{N}_d$ is non-empty. Proposition \ref{divisor-phi-circulant-prop} says that for every $n\in\mathcal{N}_d$, there exists a circulant graph of order $n$ with algebraic degree $d$. Also conversely, every $n\in\mathbb{N}$ which admits a circulant graph of the same order with algerbaic degree $d$ must belong to $\mathcal{N}_d$. With this notation, we now prove our main result regarding the order of minimal $d$-integral circulant graphs.

\begin{theorem}\label{C(d)}
Let $d\in\mathbb{N}$ with $d>1$. Then $C(d)=min(\mathcal{N}_d)$.
\end{theorem}
\pf 
For $d>1$, let $m_d=min(\mathcal{N}_d)$. Then by Proposition \ref{divisor-phi-circulant-prop}, there exists a circulant graph $\Gamma=Cay(\mathbb{Z}_{m_d},S)$ such that $Deg(\Gamma)=d$, where $|S|=\varphi(m_d)/d$.\qed 

Table \ref{tab:my_label} (obtained using SageMath \cite{stein2007sage}) lists the values of $C(d)$ and the upper bound $p_d$ in Corollary \ref{cor-circ-min-ord} for $1\le d\le 100$. Note that for $d=1$, the corresponding minimal circulant graph is $\mathcal{K}_1$. In many cases, this upper bound coincides with the actual value of $C(d)$. However, in certain cases (shown in \textcolor{red}{red} font), the value given by Theorem \ref{C(d)} is less than the bound given in Corollary \ref{cor-circ-min-ord}.




\begin{table}[H]

    \centering
\scalebox{0.9}{    \begin{tabular}{|c|c|c||c|c|c||c|c|c||c|c|c|}
    \hline
     $d$ & $C(d)$ & $p_d$ & $d$ & $C(d)$ & $p_d$ & $d$ & $C(d)$ & $p_d$ & $d$ & $C(d)$ & $p_d$ \\ \hline\hline
     $1$ & \textcolor{red}{$1$} & $3$ & $26$ & $53$ & $53$ & $51$ & $103$ & $103$ & $76$ & $457$ & $457$\\ \hline
     $2$ & $5$ & $5$ & $27$ & \textcolor{red}{$81$} & $109$ & $52$ & \textcolor{red}{$159$} & $313$ & $77$ & $463$ & $463$\\ \hline
     $3$ & $7$ & $7$ & $28$ & \textcolor{red}{$87$} & $113$ & $53$ & $107$ & $107$ & $78$ & $157$ & $157$\\ \hline
     $4$ & \textcolor{red}{$15$}& $17$ & $29$ & $59$ & $59$ & $54$ & $109$ & $109$ & $79$ & $317$ & $317$\\ \hline
     $5$ & $11$& $11$ & $30$ & $61$ & $61$ & $55$ & \textcolor{red}{$121$} & $331$ & $80$ & \textcolor{red}{$187$} & $641$\\ \hline
     $6$ & $13$& $13$ & $31$ & $311$ & $311$ & $56$ & $113$ & $113$ & $81$ & $163$ & $163$\\ \hline
     $7$ & $29$& $29$ & $32$ & \textcolor{red}{$85$} & $193$ & $57$ & $229$ & $229$ & $82$ & \textcolor{red}{$249$} & $821$\\ \hline
     $8$ & $17$& $17$ & $33$ & $67$ & $67$ & $58$ & \textcolor{red}{$177$} & $233$ & $83$ & $167$ & $167$\\ \hline
     $9$ & $19$& $19$ & $34$ & $137$ & $137$ & $59$ & $709$ & $709$ & $84$ & \textcolor{red}{$203$} & $337$\\ \hline
     $10$ & \textcolor{red}{$25$}& $41$ & $35$ & $71$ & $71$ & $60$ & \textcolor{red}{$143$} & $241$ & $85$ & $1021$ & $1021$\\ \hline
     $11$ & $23$& $23$ & $36$ & $73$ & $73$ & $61$ & $367$ & $367$ & $86$ & $173$ & $173$\\ \hline
     $12$ & \textcolor{red}{$35$} & $73$ & $37$ & $149$ & $149$ & $62$ & $373$ & $373$ & $87$ & $349$ & $349$\\ \hline
     $13$ & $53$& $53$ & $38$ & $229$ & $229$ & $63$ & $127$ & $127$ & $88$ & \textcolor{red}{$267$} & $353$\\ \hline
     $14$ & $29$& $29$ & $39$ & $79$ & $79$ & $64$ & \textcolor{red}{$255$} & $257$ & $89$ & $179$ & $179$\\ \hline
     $15$ & $31$& $31$ & $40$ & \textcolor{red}{$123$} & $241$ & $65$ & $131$ & $131$ & $90$ & $181$ & $181$\\ \hline
     $16$ & \textcolor{red}{$51$}& $97$ & $41$ & $83$ & $83$ & $66$ & \textcolor{red}{$161$} & $397$ & $91$ & $547$ & $547$\\ \hline
     $17$ & $103$& $103$ & $42$ & \textcolor{red}{$129$} & $337$ & $67$ & $269$ & $269$ & $92$ & \textcolor{red}{$235$} & $1289$\\ \hline
     $18$ & $37$& $37$ & $43$ & $173$ & $173$ & $68$ & $137$ & $137$ & $93$ & $373$ & $373$\\ \hline
     $19$ & $191$& $191$ & $44$ & $89$ & $89$ & $69$ & $139$ & $139$ & $94$ & \textcolor{red}{$849$} & $941$\\ \hline
     $20$ & $41$& $41$ & $45$ & $181$ & $181$ & $70$ & \textcolor{red}{$213$} & $281$ & $95$ & $191$ & $191$\\ \hline
     $21$ & $43$ & $43$ & $46$ & \textcolor{red}{$141$} & $277$ & $71$ & $569$ & $569$ & $96$ & $193$ & $193$\\ \hline
     $22$ & \textcolor{red}{$69$} & $89$ & $47$ & $283$ & $283$ & $72$ & \textcolor{red}{$185$} & $433$ & $97$ & $389$ & $389$\\ \hline
     $23$ & $47$ & $47$ & $48$ & $97$ & $97$ & $73$ & $293$ & $293$ & $98$ & $197$ & $197$\\ \hline
     $24$ & \textcolor{red}{$65$} & $97$ & $49$ & $197$ & $197$ & $74$ & $149$ & $149$ & $99$ & $199$ & $199$\\ \hline
     $25$ & $101$ & $101$ & $50$ & $101$ & $101$ & $75$ & $151$ & $151$ & $100$ & \textcolor{red}{$275$} & $401$\\ \hline
    \end{tabular}
    }
    \caption{Table for $C(d)$}
    \label{tab:my_label}
    
\end{table}

\section{Number of isomorphism classes of connected $d$-integral circulant graphs}\label{count-number-conn-d-int-circ-graphs}

The study of integral circulant graphs has been extensive for years. In \cite{so2006integral}, So has provided a necessary and sufficient condition for a circulant graph to be integral. Let $G_{n}(d):=\lbrace{x\in\mathbb{Z}_n:gcd(x,n)=d}\rbrace$, where $d$ is a proper divisor of $n$ (i.e., $d\neq n$) and $D(n):=\lbrace{d\in\mathbb{N}:d\mid n}\rbrace$. Then the result states that,
\begin{theorem}\label{integral-circ-graphs}[\cite{so2006integral},Theorem 7.1]
Let $\Gamma=Cay(\mathbb{Z}_n,S)$ be a circulant graph. Then $\Gamma$ is integral if and only if $S$ is the union of $G_{n}(d)$-sets, i.e., $$S=\bigcup_{d\in\mathcal{D}}G_{n}(d)$$ for some $\mathcal{D}\subseteq D(n)\setminus\lbrace{n}\rbrace$.
\end{theorem}
The $G_{n}(d)$-sets are called \textit{basic integral symbols} \cite{monius2023many}. The connection set $S$ is called an \textit{integral symbol} if $\Gamma_{S}=Cay(\mathbb{Z}_n,S)$ is integral. Thus $S$ is an integral symbol if and only if $S$ is the union of the basic integral symbols. There are a total of $\tau(n)-1$ basic integral symbols and these are used to generate all $2^{\tau(n)-1}$ integral symbols, where $\tau(\cdot)$ denotes the number of divisors of a positive integer. Thus for $n\in\mathbb{N}$, there is a one-one correspondence between the integral symbols on circulant graphs of order $n$ and the subsets of $D(n)\setminus\lbrace{n}\rbrace$. In \cite{so2006integral}, the author conjectured the following.
\begin{conjecture}
Let $\Gamma_{S}=Cay(\mathbb{Z}_n,S)$ and $\Gamma_{T}=Cay(\mathbb{Z}_n,T)$ be two integral circulant graphs. If $S\neq T$, then $sp(\Gamma_{S})\neq sp(\Gamma_{T})$, hence $\Gamma_{S}$ and $\Gamma_{T}$ are not isomorphic.
\end{conjecture}
If the conjecture is true, then it follows that \textit{there are exactly $2^{\tau(n)-1}$ non-isospectral integral circulant graphs of order $n$}. The conjecture still remains open except for some special choices of $n$, which has been resolved by M{\"o}nius and So in \cite{monius2023many}. However, the weaker version of the conjecture which states that \textit{there are exactly $2^{\tau(n)-1}$ non-isomorphic integral circulant graphs of order $n$} was proved by Klin and Kov{\'a}cs in \cite{klin2012automorphism} (Corollary 11.3). By the techniques from Schur ring and group theory, they proved that different integral symbols correspond to non-isomorphic circulant graphs and hence the number of distinct integral symbols on circulant graphs of order $n$ is equal to the number of non-isomorphic integral circulant graphs of order $n$. In fact the following result can be found in \cite{monius2021algebraic}.
\begin{theorem}\label{isomorphic-integral-circ-graphs}[\cite{monius2021algebraic},Theorem 3.3.4]
Let $\Gamma_{S}=Cay(\mathbb{Z}_n,S)$ and $\Gamma_{T}=Cay(\mathbb{Z}_n,T)$ be two integral circulant graphs. If $\Gamma_{S}\cong\Gamma_{T}$, then $S=T$.
\end{theorem}

For $d\in\mathbb{N}$, recall that $\mathcal{N}_d=\lbrace{n\in\mathbb{N}:2d\mid\varphi(n)}\rbrace$. Then by Proposition \ref{divisor-phi-circulant-prop}, for every $n\in\mathcal{N}_d$, there exists $S\subseteq\mathbb{Z}_{n}\setminus\lbrace{0}\rbrace$ with $S=-S$ such that $Deg(Cay(\mathbb{Z}_n,S))=d$. In view of this, for $d\in\mathbb{N}$, we introduce the notation $\mathcal{F}_{n,d}$ which denotes \textit{the collection of isomorphism classes of connected, $d$-integral circulant graphs of order $n\in\mathcal{N}_d$}. Note that for every $n\in\mathbb{N}$, $\mathcal{K}_n\in\mathcal{F}_{n,1}$.

In this section we find a closed form of $|\mathcal{F}_{n,1}|$ for all $n\in\mathbb{N}$ and provide some bounds on $|\mathcal{F}_{n,d}|$ for $d>1$ and $n\in\mathcal{N}_d$. Also we find an explicit formula for $|\mathcal{F}_{p,d}|$ where $d$ is a prime and $p$ is a prime such that $p\in\mathcal{N}_d$.

\subsection{\textbf{Closed form of $\mathbf{|\mathcal{F}_{n,1}|}$}}


\begin{proposition}\label{F_{n,1}}
For $n>1$, $|\mathcal{F}_{n,1}|\ge 2^{\tau(n)-2}$. The equality holds if and only if $n$ is a prime power.
\end{proposition}
\pf Note that $G_{n}(1)=\mathbb{Z}_{n}^{*}$. Hence if the connection set $S$ contains $G_{n}(1)$, the graph $Cay(\mathbb{Z}_n,S)$ must be connected. Since there are $\tau(n)-2$ many $G_{n}(d)$ sets, where $d\in D(n)\setminus\lbrace{1,n}\rbrace$, we can choose $k$ of them in $\tau(n)-2\choose k$ many ways for $k=0,1,\ldots,\tau(n)-2$. Hence the total number of possible ways to choose the $G_{n}(d)$-sets is $2^{\tau(n)-2}$, where $d\in D(n)\setminus\lbrace{1,n}\rbrace$. Thus the circulant graph $Cay(\mathbb{Z}_n,S)$ becomes connected, where $$S=G_{n}(1)\cup\left(\bigcup_{d\in\mathcal{D}}G_{n}(d)\right)$$ for some $\mathcal{D}\subseteq D(n)\setminus\lbrace{1,n}\rbrace$. By the means of Theorem \ref{integral-circ-graphs} and Theroem \ref{isomorphic-integral-circ-graphs}, the result follows.

If $n=p^r$, the proper divisors of $n$ are $1,p,p^2,\ldots,p^{r-1}$. Since $\bigcup_{d\in\mathcal{D}}G_{n}(d)$ where $\mathcal{D}\subseteq D(n)\setminus\lbrace{1,n}\rbrace$, gives rise to a connection set $S$ for which the corresponding circulant graph $Cay(\mathbb{Z}_n,S)$ is disconnected, we ensure by our previous argument that there are exactly $2^{\tau(n)-2}$ connected, non-isomorphic integral circulant graphs of order $n=p^r$.

Conversely, let $|\mathcal{F}_{n,1}|=2^{\tau(n)-2}$. Suppose $n$ has two distinct prime divisors, say $p$ and $q$. Consider the connection set $T:=G_{n}(p)\cup G_{n}(q)$. Then clearly $\langle{T}\rangle=\mathbb{Z}_n$. Then $Cay(\mathbb{Z}_n,T)$ is a connected, integral circulant graph, where the connection set $T$ does not contain $G_{n}(1)$. Thus there must be at least $2^{\tau(n)-2}+1$ non-isomorphic, integral circulant graphs of order $n$, which contradicts our initial assumption. Hence $n$ must be some prime power.\qed

Before going into the next result, we first introduce some notations. 
For $n\in\mathbb{N}$, let $\mathcal{S}(n)$ be the collection of all integral symbols on circulant graphs of order $n$. Clearly, $|\mathcal{S}(n)|=2^{\tau(n)-1}$. For each divisor $d$ of $n$, let $\mathcal{S}_{C}(d)$ denote the integral symbol on connected circulant graphs of order $d$. Define $\mathcal{T}(n):=\bigcup_{d\mid n}\mathcal{S}_{C}(d)$. For $A\subseteq\mathbb{N}$, define $n/A:=\lbrace{\frac{n}{a}:a\in A}\rbrace$. Let $[A]$ denote the lcm of the integers in $A$ and $(A)$ denote the gcd of the integers in $A$.

\begin{lemma}\label{gcd-lemma}
Let $\emptyset\neq A\subseteq\mathbb{N}$ be finite. Then $([A]/A)=1$.
\end{lemma}
\pf Let $A:=\lbrace{a_1,\ldots,a_k}\rbrace$. Let $b_i:=\Pi_{j\neq i}a_j$, for $i=1,\ldots,k$ and $M:=gcd(b_1,\ldots,b_k)$. Then $[A]=lcm(a_1,\ldots,a_k)=\frac{a_1\cdots a_k}{M}$. Hence
$$([A]/A)=gcd\left(\frac{[A]}{a_1},\ldots,\frac{[A]}{a_k}\right)=gcd\left(\frac{b_1}{M},\ldots,\frac{b_k}{M}\right)=1.$$\qed

\begin{lemma}\label{Lemma-S(n)-T(n)}
For every $n\in\mathbb{N}$, the sets $\mathcal{S}(n)$ and $\mathcal{T}(n)$ have the same size.
\end{lemma}
\pf Define the function $\Phi:\mathcal{S}(n)\to\mathcal{T}(n)$ in the following way.
$$\Phi:\bigcup_{d\in\mathcal{D}}G_{n}(d)\mapsto\bigcup_{a\in n/\mathcal{D}}G_{[n/\mathcal{D}]}\left(\frac{[n/\mathcal{D}]}{a}\right)$$
where $\mathcal{D}\subseteq D(n)\setminus\lbrace{n}\rbrace$. If $\mathcal{D}=\emptyset$, then $\Phi$ maps the empty symbol to the empty symbol i.e., the corresponding circulant graph is an empty graph. 
More precisely, $\Phi$ maps the empty circulant graph on $n$ vertices (i.e., $\bigcup_{i=1}^n\mathcal{K}_1=n\mathcal{K}_1$, i.e., $n$ many isolated points) to the empty circulant graph on single vertex (i.e., $\mathcal{K}_1$).
For $\emptyset\neq\mathcal{D}\subseteq D(n)\setminus\lbrace{n}\rbrace$, let $\mathcal{H}:=[A]/A$, where $A=n/\mathcal{D}$. By Lemma \ref{gcd-lemma}, $(\mathcal{H})=1$. Hence $\langle{S}\rangle=\mathbb{Z}_{[n/\mathcal{D}]}$, where $$S:=\bigcup_{a\in n/\mathcal{D}}G_{[n/\mathcal{D}]}\left(\frac{[n/\mathcal{D}]}{a}\right).$$
Thus $\Phi\left(\bigcup_{d\in\mathcal{D}}G_{n}(d)\right)$ becomes an integral symbol on connected circulant graphs of order $[n/\mathcal{D}]$. Hence $\Phi$ is well-defined. We prove that $\Phi$ is a bijection.

Let $\Phi(\mathcal{A}_1)=\Phi(\mathcal{A}_2)$, for some $\mathcal{A}_1,\mathcal{A}_2\in\mathcal{S}(n)$. Since $\mathcal{A}_1,\mathcal{A}_2$ are integral symbols, let $\mathcal{A}_1=\bigcup_{d\in\mathcal{D}_1}G_{n}(d)$ and $\mathcal{A}_2=\bigcup_{d\in\mathcal{D}_2}G_{n}(d)$, for some $\mathcal{D}_1,\mathcal{D}_2\subseteq D(n)\setminus\lbrace{n}\rbrace$. Hence
$$\bigcup_{a\in n/\mathcal{D}_1}G_{[n/\mathcal{D}_1]}\left(\frac{[n/\mathcal{D}_1]}{a}\right)=\bigcup_{a\in n/\mathcal{D}_2}G_{[n/\mathcal{D}_2]}\left(\frac{[n/\mathcal{D}_2]}{a}\right).$$
Since the left-hand side and the right-hand side of the equality are subsets of $\mathbb{Z}_{[n/\mathcal{D}_1]}$ and $\mathbb{Z}_{[n/\mathcal{D}_2]}$ respectively, it follows that $[n/\mathcal{D}_1]=[n/\mathcal{D}_2]$. For brevity, we write $[n/\mathcal{D}_i]=M$ for $i=1,2$. Let $k\in n/\mathcal{D}_1$. 
Since the $G_{n}(d)$-sets are disjoint, $G_{M}\left(\frac{M}{k}\right)=G_{M}\left(\frac{M}{l}\right)$, for some $l\in n/\mathcal{D}_2$. Thus $\frac{M}{k}\in G_{M}\left(\frac{M}{l}\right)$, i.e., $gcd\left(\frac{M}{k},M\right)=\frac{M}{l}$ and hence $k=l$. Therefore, $n/\mathcal{D}_1\subseteq n/\mathcal{D}_2$. Similarly the reverse inclusion can also be shown. Thus $n/\mathcal{D}_1=n/\mathcal{D}_2$, i.e., $\mathcal{D}_1=\mathcal{D}_2$. Hence $\mathcal{A}_1=\mathcal{A}_2$ and $\Phi$ is injective.

Let $\emptyset\neq\mathcal{A}\in\mathcal{T}(n)$. Then $\mathcal{A}\in\mathcal{S}_{C}(t)$, for some $t\mid n$. Therefore $\mathcal{A}$ is of the form $\bigcup_{s\in\mathcal{D}}G_{t}(s)$, for some $\mathcal{D}\subseteq D(t)\setminus\lbrace{t}\rbrace$. Let $\mathcal{B}:=t/\mathcal{D}$. Since $\mathcal{A}$ is a connected integral symbol, $(\mathcal{D})=1$. A little calculation shows that $[\mathcal{B}]=[t/\mathcal{D}]=t$. Let $$\mathcal{M}:=\bigcup_{d\in n/\mathcal{B}}G_{n}(d).$$ Then $\mathcal{M}\in\mathcal{S}(n)$, as $n/\mathcal{B}\subseteq D(n)\setminus\lbrace{n}\rbrace$. Thus we have $$\Phi(\mathcal{M})=\bigcup_{a\in\mathcal{B}}G_{[\mathcal{B}]}\left(\frac{[\mathcal{B}]}{a}\right)=\bigcup_{a\in t/\mathcal{D}}G_{t}\left(\frac{t}{a}\right)=\bigcup_{s\in\mathcal{D}}G_{t}(s)=\mathcal{A}.$$ Hence $\Phi$ is surjective.\qed

\begin{theorem}\label{Thm-F{n,1}}
For $n\in\mathbb{N}$,
$$|\mathcal{F}_{n,1}|=\frac{1}{2}\sum_{d\mid n}\mu(d)2^{\tau(n/d)}$$
where $\mu(\cdot)$ is the M{\"o}bius function.
\end{theorem}
\pf By Theorem \ref{isomorphic-integral-circ-graphs}, we have $|\mathcal{S}_{C}(d)|=|\mathcal{F}_{d,1}|$ for $d\mid n$. Also from Lemma \ref{Lemma-S(n)-T(n)}, we have $|\mathcal{S}(n)|=|\mathcal{T}(n)|$.  Since the sets $\mathcal{S}_{C}(d)$ are mutually disjoint for distinct divisors $d$ of $n$, we have $|\mathcal{T}(n)|=\sum_{d\mid n}|\mathcal{S}_{C}(d)|$ i.e., $\sum_{d\mid n}|\mathcal{F}_{d,1}|=2^{\tau(n)-1}$. By M{\"o}bius inversion formula, we obtain our desired result.\qed

Theorem \ref{Thm-F{n,1}} leads us to the following immediate corollary, which is a known result from \cite{monius2020algebraic} (Corollary 2.6).

\begin{corollary}\label{integral-circ-graph-prime-ord}
For a prime $p$, $|\mathcal{F}_{p,1}|=1$. In other words, the only connected, integral circulant graph on $p$ vertices is $\mathcal{K}_p$. The only other (disconnected) integral circulant graph on $p$ vertices is the empty graph $\bigcup_{i=1}^p\mathcal{K}_1=p\mathcal{K}_1$.
\end{corollary}

\subsection{\textbf{Bounds on $\mathbf{|\mathcal{F}_{n,d}|~for~d>1~and~n\in\mathcal{N}_d}$}}
To begin with the bound of $|\mathcal{F}_{n,d}|$ where $d>1$ and $n\in\mathcal{N}_d$, we need the following lemmas first. Let $G$ be a group and $X\subseteq G$. We say that an element $a\in G$ fixes $X$ if $aX=X$ where $aX=\lbrace{ax:x\in X}\rbrace$.

\begin{lemma}\label{group-subset-lemma}
Let $G$ be a finite group and let $X\subseteq G$ be such that $e\in X$ but $X$ is not a subgroup of $G$. If $X$ has $p$ elements, for some prime $p$, then $e$ is the only element in $G$ that fixes $X$. In other words, if $aX=X$ for some $a\in G$, then $a=e$.
\end{lemma}

\pf Define $S:=\lbrace{a\in G:aX=X}\rbrace$. Then $S$ is a subgroup of $G$. Since $e\in X$, no elements outside of $X$ fix $X$. Hence $S\subseteq X$. Let $\psi:S\times X\to X$ be the left multiplication action defined as $\psi(y,x)=yx$. By Burnside lemma, $$\#Orbits(\psi)=\frac{1}{|S|}\sum_{y\in S}|Fix(y)|$$ where $Fix(y)=\lbrace{x\in X:\psi(y,x)=x}\rbrace$. Since $Fix(e)=X$ and $Fix(y)=\emptyset$, for $y\neq e$, we have $$\#Orbits(\psi)=\frac{1}{|S|}\left(|Fix(e)|+\sum_{\tiny\begin{array}{cc}
    y\in S\\
    y\neq e 
\end{array}}|Fix(y)|\right)=\frac{|X|}{|S|}.$$ Since $|X|$ is prime, $|S|=1$ or $|X|$. But the latter cannot hold, since $X$ is not a subgroup of $G$. Hence $S=\lbrace{e}\rbrace$.\qed

\begin{lemma}\label{group-subset-lemma-2}
Let $G$ be a finite group and let $X\subseteq G$. If $gcd(|X|,|G|)=1$, then $e$ is the only element in $G$ that fixes $X$.
\end{lemma}
\pf From Lemma \ref{group-subset-lemma}, $|S|$ divides both $|X|$ and $|G|$. As $gcd(|X|,|G|)=1$, $S=\lbrace{e}\rbrace$.\qed

\begin{theorem}\label{counting-d-int-circ-graph}
Let $d>1$ and  $n\in\mathcal{N}_d$. Then\\
$$|\mathcal{F}_{n,d}|\ge\begin{cases}
    \varphi(d), \mbox{ if d is a prime power. }\\
    \\
    \varphi(d)+\omega(d), \mbox{ if $d$ is not a prime power. }\\
\end{cases}$$
where $\omega(\cdot)$ denotes the number of distinct prime divisors of a positive integer.
\end{theorem}
\pf Let $d>1$ and $n\in\mathcal{N}_d$. By Proposition \ref{divisor-phi-circulant-prop}, there exists a $d$-integral circulant graph $\Gamma=Cay(\mathbb{Z}_n,H)$ such that $H\le\mathbb{Z}_n^{*}$ with $|H|=\frac{\varphi(n)}{d}$ and $H$ is inverse-symmetric in $\mathbb{Z}_n$. We define the following subsets of $\mathbb{Z}_n^{*}/H$ inductively.
$$\mathcal{T}_1:=\lbrace{H}\rbrace$$
$$\mathcal{T}_2:=\lbrace{H,x_2H}\rbrace$$
$$\mathcal{T}_3:=\lbrace{H,x_2H,x_3H}\rbrace$$
$$\vdots$$
$$\mathcal{T}_{d-1}:=\lbrace{H,x_2H,x_3H,\ldots,x_{d-1}H}\rbrace$$
$$\mathcal{T}_d:=\lbrace{H,x_2H,x_3H,\ldots,x_dH}\rbrace$$
for some $x_2,x_3,\ldots,x_d\in\mathbb{Z}_{n}^{*}$, the cosets being mutually disjoint. Define $S_1:=H$ and
$$S_m:=H\cup\left(\bigcup_{i=2}^{m}x_iH\right)$$
for $m=2,\ldots,d$. Clearly for each $1\le m\le d$, $S_m$ is inverse-symmetric in $\mathbb{Z}_n$ and $|S_m|=m|H|=\frac{m\varphi(n)}{d}$. Since $\mathcal{T}_d=\mathbb{Z}_n^{*}/H$, $S_d=\mathbb{Z}_{n}^{*}$. For each $1\le m\le d$, define the circulant graphs $\Gamma_{S_m}:=Cay(\mathbb{Z}_n,S_m)$. For $m=d$, the circulant graph $\Gamma_{S_d}$ becomes $Cay(\mathbb{Z}_n,\mathbb{Z}_{n}^{*})$ which is integral by \ref{formula-alg-deg-circ}. Therefore, we discard the case $m=d$. From the construction of $S_m$, it is clear that $S_m$ is fixed by $H$ and the circulant graphs $\Gamma_{S_m}$ are connected and non-isomorphic with valency $\frac{m\varphi(n)}{d}$ for each $1\le m\le d-1$. Define $Fix(S_m):=\lbrace{a\in\mathbb{Z}_{n}^{*}:aS_m=S_m}\rbrace$. Then $Fix(S_m)$ is the maximum subgroup of $\mathbb{Z}_{n}^{*}$ fixing $S_m$. Clearly, $H\subseteq Fix(S_m)$. The rest of the proof is divided into two cases.\\
\textbf{Case 1:} $d$ is a prime power.

Take $1\le m\le d-1$ such that $gcd(m,d)=1$. Since $\mathcal{T}_m\subseteq\mathbb{Z}_{n}^{*}/H$ and $|\mathcal{T}_m|=m$, by Lemma \ref{group-subset-lemma-2}, $H$ is the only element that fixes $\mathcal{T}_m$. Hence for each $1\le m\le d-1$ with $gcd(m,d)=1$, $Fix(S_m)=H$ and thus by \ref{formula-alg-deg-circ}, $Deg(\Gamma_{S_m})=d$. Hence $$|\mathcal{F}_{n,d}|\ge\varphi(d).$$
\textbf{Case 2:} $d$ is not a prime power.

Then there exist primes $p_1,p_2$ such that $2\le p_1<p_2$ and $p_1,p_2\mid d$. By Cauchy's theorem, in the above construction of $\mathcal{T}_m$'s, choose elements $x_2H,x_3H\in\mathbb{Z}_{n}^{*}/H$ such that $o(x_2H)=p_2$ and $o(x_3H)=p_1$. Clearly, $\mathcal{T}_2=\lbrace{H,x_2H}\rbrace$ and $\mathcal{T}_3=\lbrace{H,x_2H,x_3H}\rbrace$ are not subgroups of $\mathbb{Z}_{n}^{*}/H$ and similarly for all primes $q\mid d$, $\mathcal{T}_q$ is not a subgroup of $\mathbb{Z}_{n}^{*}/H$. Since $H\in\mathcal{T}_q\subseteq\mathbb{Z}_{n}^{*}/H$, by Lemma \ref{group-subset-lemma}, $H$ is the only element that fixes $\mathcal{T}_q$ and thus $Fix(S_q)=H$ for all primes $q\mid d$. Also by Lemma \ref{group-subset-lemma-2}, $Fix(S_m)=H$ for each $1\le m\le d-1$ with $gcd(m,d)=1$.
Thus by \ref{formula-alg-deg-circ}, we have $Deg(\Gamma_{S_q})=Deg(\Gamma_{S_m})=d$, for all primes $q\mid d$ and for all $1\le m\le d-1$ with $gcd(m,d)=1$. Hence $$|\mathcal{F}_{n,d}|\ge\varphi(d)+\omega(d).$$\qed

From Theorem \ref{counting-d-int-circ-graph}, we have, for a positive integer $d>2$ and $n\in\mathcal{N}_d$, $|\mathcal{F}_{n,d}|\ge 2$. Thus, in particular, for $d>2$, $|\mathcal{F}_{C(d),d}|\ge 2$. This leads us to the following immediate corollary.

\begin{corollary}
For $d>2$, the minimal $d$-integral circulant graph is not unique. 
\end{corollary}

From Table \ref{tab:my_label}, $C(1)=1$ and $C(2)=5$. Thus from Theorem \ref{Thm-F{n,1}} and from Theorem \ref{exact-value-F_{p,d}}, it follows respectively that $|\mathcal{F}_{1,1}|=1$ and $|\mathcal{F}_{5,2}|=1$. Thus for $d=1$ and $2$, the minimal $d$-integral circulant graphs are unique (up to isomorphism). In fact it is easily noticeable that $\mathcal{K}_1$ is the unique minimal integral circulant graph and $\mathcal{C}_5$ is the unique minimal $2$-integral circulant graph.

\begin{remark}
The construction of the circulant graphs $\Gamma_{S_m}$ in Theorem \ref{counting-d-int-circ-graph} suggests that there are at least $\varphi(d)$ many non-isospectral, connected, $d$-integral circulant graphs of order $n\in\mathcal{N}_d$. 
\end{remark}

The lower bound in Theorem \ref{counting-d-int-circ-graph} can be enhanced if we consider $d$ to be square free, which we show in our next theorem.

\begin{theorem}\label{sq-free-LB}
Let $d>1$ be a square free integer and $n\in\mathcal{N}_d$. Then $|\mathcal{F}_{n,d}|\ge d-\omega(d)$.
\end{theorem}

\pf We adapt the same notations from Theorem \ref{counting-d-int-circ-graph} for convenience. Let $d>1$ be a square free integer and $n\in\mathcal{N}_d$. If $d$ is prime, then $\omega(d)=1$. Thus $|\mathcal{F}_{n,d}|\ge d-1$, which is already true by Theorem \ref{counting-d-int-circ-graph}. Hence we consider $d$ to be composite.

Let $\omega(d)=k\ge 2$. Then $d=p_1p_2\cdots p_k$, where $p_i$'s are primes with $p_1<p_2<\cdots<p_k$. By Lemma \ref{divisor-phi-circulant-lemma}, let $H\le\mathbb{Z}_{n}^{*}$ be inverse-symmetric in $\mathbb{Z}_n$ and $|H|=\frac{\varphi(n)}{d}$. Since $d$ is square free, $\mathbb{Z}_{n}^{*}/H$ is cyclic of order $d$. Thus for every $t\mid d$, there exist exactly $\varphi(t)$ elements of order $t$ in $\mathbb{Z}_{n}^{*}/H$. We now construct the subsets of $\mathbb{Z}_{n}^{*}/H$ inductively.

Define $\mathcal{T}_1:=\lbrace{H}\rbrace$ and then construct the subsets $\mathcal{T}_m$ by adjoining only $(\varphi(p_i)-1)$-many elements of $\mathbb{Z}_{n}^{*}/H$ of order $p_i$ to $\mathcal{T}_1$ one after another, for each $i=1,\ldots,k$. Once all the elements of prime orders have been used, we start adjoining the elements of non-prime order $t\mid d$ with $t\neq 1$, one after another. Unlike the earlier case, here we consider all the $\varphi(t)$ elements of order $t$. Thus in the entire process, all the elements of order $t\mid d$ have been exhausted in constructing the subsets $\mathcal{T}_m$, except one element of order $p_i$ for each $i=1,\ldots,k$. Hence the entire process yields exactly $d-k$ subsets of $\mathbb{Z}_{n}^{*}/H$.

Observe that the way the subsets $\mathcal{T}_m$ of $\mathbb{Z}_{n}^{*}/H$ are constructed, for $1\le m\le d-k$, no $\mathcal{T}_m$ contains any subgroup of $\mathbb{Z}_{n}^{*}/H$ of order $p_i$. Thus for $1\le m\le d-k$, no $\mathcal{T}_m$ contains any subgroup of $\mathbb{Z}_{n}^{*}/H$ of order $t\mid d$, except the identity subgroup $\lbrace{H}\rbrace$. Therefore, for each $1\le m\le d-k$, we can conclude that $H$ is the only element that fixes $\mathcal{T}_m$. Now construct the circulant graphs $\Gamma_{S_m}:=Cay(\mathbb{Z}_n,S_m)$ similarly as in Theorem \ref{counting-d-int-circ-graph}, where $$S_m:=\bigcup_{\mathcal{A}\in\mathcal{T}_m}\mathcal{A}.$$ Thus $Fix(S_m)=H$ and hence by \ref{formula-alg-deg-circ}, $Deg(\Gamma_{S_m})=d$ for $m=1,2,\ldots,d-k$. Hence $$|\mathcal{F}_{n,d}|\ge d-\omega(d).$$\qed

\subsection{\textbf{Bounds on $\mathbf{|\mathcal{F}_{p,d}|~for~d>1~and~p\in\mathcal{N}_d}$}}

In this section, we consider circulant graphs of prime order. We omit using the term `connected' while considering $d$-integral circulant graphs of order $p$, since every non-empty circulant graph of prime order is essentially connected.
We start by determining a lower bound of $|\mathcal{F}_{p,d}|$.

\begin{theorem}\label{counting-prime-ord-d-alg-deg-non-iso-circ-graphs}
Let $p>3$ be a prime. Then $|\mathcal{F}_{p,d}|\ge d-1$, for all $d\mid \frac{p-1}{2}$ with $d\neq 1$.
\end{theorem}
\pf The main idea of this proof lies in Theorem \ref{counting-d-int-circ-graph}. Let $p>3$ be a prime. For $d\mid\frac{p-1}{2}$ with $d\neq 1$, let $H$ be the unique subgroup of $\mathbb{Z}_{p}^{*}$ of order $\frac{p-1}{d}$. Then $Deg(Cay(\mathbb{Z}_p,H))=\frac{\varphi(p)}{|H|}=d$ by \ref{formula-alg-deg-circ}. 
Note that $H=\lbrace{x\in\mathbb{Z}_{p}^{*}:x^{(p-1)/d}\equiv 1~(mod~p)}\rbrace$. 
If $r$ is a primitive root in $\mathbb{Z}_{p}^{*}$, then $H=\langle{r^d}\rangle$.  Thus $\mathbb{Z}_{p}^{*}/H=\lbrace{H,rH,\ldots,r^{d-1}H}\rbrace$. Define $$S_m:=\bigcup_{i=0}^{m-1}r^{i}H$$ for $m=1,\ldots,d-1$. Then $S_m$ is inverse-symmetric in $\mathbb{Z}_p$ and $|S_m|=m|H|=\frac{m(p-1)}{d}$ for $1\le m\le d-1$. Consider the circulant graphs $\Gamma_{S_m}:=Cay(\mathbb{Z}_p,S_m)$. By the construction of $S_m$, $H$ fixes $S_m$ for each $m$ and the circulant graphs $\Gamma_{S_m}$ are non-isomorphic with valency $\frac{m(p-1)}{d}$. 

For $1\le m\le d-1$, consider the following sum,
$$A_m=\sum_{j\in S_m}j^{(p-1)/d}=\sum_{h\in H}h^{(p-1)/d}+r^{(p-1)/d}\sum_{h\in H}h^{(p-1)/d}+
\cdots+r^{(m-1)(p-1)/d}\sum_{h\in H}h^{(p-1)/d}$$
$$=\left(1+r^{(p-1)/d}+
\cdots+r^{(m-1)(p-1)/d}\right)\sum_{h\in H}h^{(p-1)/d}$$
$$~~=\left[\frac{r^{m(p-1)/d}-1}{r^{(p-1)/d}-1}\right]|H|=\left[\frac{r^{m(p-1)/d}-1}{r^{(p-1)/d}-1}\right]\left(\frac{p-1}{d}\right).$$
Recall $Fix(S_m)$ from Theorem \ref{counting-d-int-circ-graph}. Let $a\in Fix(S_m)$.
Then
$$a^{(p-1)/d}A_m\equiv a^{(p-1)/d}\sum_{j\in S_m}j^{(p-1)/d}\equiv \sum_{j\in S_m}(aj)^{(p-1)/d}\equiv \sum_{l\in S_m}l^{(p-1)/d}\equiv A_m~(mod~p).$$
Thus $(a^{(p-1)/d}-1)A_m\equiv 0~(mod~p)$. Since $r$ is a primitive root, it can be shown that $p\nmid A_m$ and hence
$a^{(p-1)/d}\equiv 1~(mod~p)$. Therefore $a\in H$ and hence $Fix(S_m)=H$, for each $1\le m\le d-1$. By \ref{formula-alg-deg-circ}, $Deg(\Gamma_{S_m})=d$. Thus there are at least $(d-1)$-many non-isomorphic circulant graphs on $p$ vertices with algebraic degree $d$, completing the proof.\qed

In \cite{monius2020algebraic}, the author showed that (Corollary 2.8) for a prime $p>2$, there exist at least $\varphi(\frac{p-1}{2})$ non-isomorphic circulant graphs on $p$ vertices and of maximum algebraic degree within the family of all circulant graphs. Therefore, according to our notation $\left|\mathcal{F}_{p,\frac{p-1}{2}}\right|\ge\varphi(\frac{p-1}{2})$. The following result which is a corollary to Theorem \ref{counting-prime-ord-d-alg-deg-non-iso-circ-graphs}, improves the lower bound by a significant margin.

\begin{corollary}
Let $p>3$ be a prime. Then
$$\left|\mathcal{F}_{p,\frac{p-1}{2}}\right|\ge \frac{p-3}{2}.$$
\end{corollary}



Next we provide an upper bound on $|\mathcal{F}_{p,d}|$. The following is a well-known result regarding the isomorphism between circulant graphs.

\begin{theorem}\label{sufficient-circ-isomorphism-theorem}
Let $\Gamma_S=Cay(\mathbb{Z}_n,S)$ and $\Gamma_T=Cay(\mathbb{Z}_n,T)$. If there exists $m\in\mathbb{Z}_{n}^{*}$ such that $S=mT$, then $\Gamma_S\cong\Gamma_T$.
\end{theorem}

\begin{remark}\label{Adam_conjecture}
    The converse of Theorem \ref{sufficient-circ-isomorphism-theorem}, i.e., whether $\Gamma_S\cong\Gamma_T$ implies $S=mT$ for some $m\in\mathbb{Z}_{n}^{*}$, is not true in general. However, P.P. Pálfy in \cite{palfy1987isomorphism}, proved that the converse holds if $gcd(n,\varphi(n))=1$. Thus, in particular, the converse holds for circulant graphs of prime order. 
\end{remark}

We state a lemma which is needed for proving the next theorem. We omit the proof considering its triviality.

\begin{lemma}\label{S-T-lemma}
Let $G$ be a finite abelian group and $H\le G$. Let $\mathcal{T}:=\lbrace{H,x_2H,x_3H,\ldots,x_nH}\rbrace$ be a subset of $G/H$ and $S:=H\cup x_2H\cup x_3H\cup\cdots\cup x_nH$ be a subset of $G$ for some $x_2,x_3,\ldots,x_n\in G$. Then $\mathcal{T}$ is a subgroup of $G/H$ if and only if $S$ is a subgroup of $G$.
\end{lemma}

\begin{theorem}\label{counting-prime-ord-d-UB}
Let $d\in\mathbb{N}$ with $d>1$. Then for any prime $p\in\mathcal{N}_d$, $$|\mathcal{F}_{p,d}|\le\frac{1}{d}(2^d-2)+\frac{1}{d}(d-\varphi(d)-1)\left(\sum_{\tiny\begin{array}{cc}
     1\le m\le d-1  \\
     gcd(m,d)>1
\end{array}}{d\choose m}\right)-(\sigma(d)-1-d)$$
where $\sigma(\cdot)$ denotes the sum of divisors of a positive integer.
\end{theorem}
\pf Let $d>1$ and $p\in\mathcal{N}_d$ be prime. Let $\Gamma_S:=Cay(\mathbb{Z}_p,S)$ be a circulant graph such that $Deg(\Gamma_S)=d$. By Theorem 6.1 in \cite{monius2022splitting}, the maximum subgroup of $\mathbb{Z}_{p}^{*}$ fixing $S$ must be of order $\frac{p-1}{d}$. From the uniqueness of a subgroup of $\mathbb{Z}_{p}^{*}$, this subgroup must be the subgroup $H$ defined in Theorem \ref{counting-prime-ord-d-alg-deg-non-iso-circ-graphs}. Recall from Theorem \ref{counting-prime-ord-d-alg-deg-non-iso-circ-graphs} that $\mathbb{Z}_{p}^{*}/H=\lbrace{H,rH,\ldots,r^{d-1}H}\rbrace$ where $r$ is a primitive root in $\mathbb{Z}_{p}^{*}$. Since $S$ is fixed by $H$, $S$ is a union of some $H$-orbits. Thus $S$ can be written as $$S=\bigcup_{j\in D}r^jH$$ for some $D\subseteq\lbrace{0,1,\ldots,d-1}\rbrace$ with $1\le |D|\le d-1$. Note that if $|D|=d$, then $|S|=p-1$ implying $Deg(\Gamma_S)=1$ which opposes our assumption of $d>1$. Let $\mathcal{R}_m$ denote the set of all $m$-subsets of $\mathbb{Z}_{p}^{*}/H$ for $1\le m\le d-1$. Define the following action on $\mathcal{R}_m$.
$$\psi_m:\mathbb{Z}_{p}^{*}/H\times\mathcal{R}_m\to\mathcal{R}_m$$
$$(r^iH,A)\mapsto(r^iH)A$$
The orbits of this action are given by 
$Orb^{(\psi_m)}(A)=\lbrace{(r^iH)A:0\le i\le d-1}\rbrace$ for $A\in\mathcal{R}_m$. If $Orb^{(\psi_m)}(A)\neq Orb^{(\psi_m)}(B)$ for some $A,B\in\mathcal{R}_m$, then $A\neq (uH)B$ for any $u\in\mathbb{Z}_{p}^{*}$ and conversely. Thus it follows from Remark \ref{Adam_conjecture} that each of the orbits of the action $\psi_m$ actually denotes each isomorphism class of circulant graphs of order $p$ with valency $\frac{m(p-1)}{d}$ having algebraic degree at most $d$.

Now $Fix^{(\psi_m)}(r^iH)=\lbrace{A\in\mathcal{R}_m:(r^iH)A=A}\rbrace\subseteq\mathcal{R}_m$. Clearly, $Fix^{(\psi_m)}(H)=\mathcal{R}_m$. We claim that $Fix^{(\psi_m)}(r^iH)=\emptyset$ for every $0\le i\le d-1$ with $gcd(i,d)=1$.

Suppose $Fix^{(\psi_m)}(r^iH)\neq\emptyset$. Then there exists $A\in\mathcal{R}_m$ such that $r^iH$ fixes $A$ i.e., $(r^iH)A=A$. It follows that the subgroup $\langle{r^iH}\rangle$ of $\mathbb{Z}_p^{*}/H$ also fixes $A$. Then by the argument in Lemma \ref{group-subset-lemma}, $|\langle{r^iH}\rangle|$ divides $|A|$. Since $gcd(i,d)=1$, $|\langle{r^iH}\rangle|=d$ which implies $|A|\ge d$, a contradiction.

Thus by Burnside Lemma, for each $1\le m\le d-1$,
$$\#Orbits(\psi_m)=\frac{1}{d}\sum_{i=0}^{d-1}|Fix^{(\psi_m)}(r^iH)|=\frac{1}{d}\left(\sum_{i\in\mathbb{Z}_d^{*}}|Fix^{(\psi_m)}(r^iH)|+\sum_{i\notin\mathbb{Z}_d^{*}}|Fix^{(\psi_m)}(r^iH)|\right)$$
$$~~~~~~~~~~~~~~~=\frac{1}{d}\sum_{i\notin\mathbb{Z}_d^{*}}|Fix^{(\psi_m)}(r^iH)|.$$
Now take $1\le m\le d-1$ such that $gcd(m,d)=1$. Then by Lemma \ref{group-subset-lemma-2}, it follows that $Fix^{(\psi_m)}(r^iH)=\emptyset$ for $i\neq 0$. 
Hence for $1\le m\le d-1$ with $gcd(m,d)=1$, we have,
$$\#Orbits(\psi_m)=\frac{1}{d}\left(|Fix^{(\psi_m)}(H)|+\sum_{\scriptsize\begin{array}{cc}
    i\notin\mathbb{Z}_d^{*} \\
    i\neq 0
\end{array}}|Fix^{(\psi_m)}(r^iH)|\right)=\frac{1}{d}{d\choose m}.$$
If $1\le m\le d-1$ be such that $gcd(m,d)>1$, then we obtain,
$$\#Orbits(\psi_m)=\frac{1}{d}\sum_{i\notin\mathbb{Z}_d^{*}}|Fix^{(\psi_m)}(r^iH)|\le\frac{1}{d}(d-\varphi(d)){d\choose m}.$$
Thus for each $1\le m\le d-1$ with $gcd(m,d)=1$, there are exactly $\frac{1}{d}{d\choose m}$ non-isomorphic circulant graphs of order $p$ with valency $\frac{m(p-1)}{d}$ having algebraic degree exactly $d$ (the assertion of the algebraic degree being {\it exactly} $d$ follows from Lemma \ref{group-subset-lemma-2} and \ref{formula-alg-deg-circ}) and for each $1\le m\le d-1$ with $gcd(m,d)>1$, there are at most $\frac{1}{d}(d-\varphi(d)){d\choose m}$ non-isomorphic circulant graphs of order $p$ with valency $\frac{m(p-1)}{d}$ having algebraic degree at most $d$. Thus total number of orbits is bounded above by the following expression,
$$\#Orbits(\psi_m)_{m=1}^{d-1}=\#Orbits(\psi_m)_{gcd(m,d)=1}+\#Orbits(\psi_m)_{gcd(m,d)>1}$$
$$~~~~~~~~~~~~~~~~~~~~~~~~~~~~~~~~~~~~~~~~~~~\le\frac{1}{d}\left(\sum_{\tiny\begin{array}{cc}
     1\le m\le d-1  \\
     gcd(m,d)=1
\end{array}}{d\choose m}\right)+\frac{1}{d}(d-\varphi(d))\left(\sum_{\tiny\begin{array}{cc}
     1\le m\le d-1  \\
     gcd(m,d)>1
\end{array}}{d\choose m}\right).$$
We now show that the expression in the second summand contains some isomorphism classes of circulant graphs of order $p$ whose algebraic degree is less than $d$, when $d$ is composite. Note that if $d$ is prime, then the second summand actually does not exist. Therefore, we choose $d$ to be composite for the rest of the proof.

Let $k\mid d$ such that $k\neq 1,d$. Let $A_k:=\langle{r^kH}\rangle\le\mathbb{Z}_{p}^{*}/H$. Then $|A_k|=d/k$ and hence $A_k\in\mathcal{R}_{d/k}$. If $$S:=\bigcup_{i=0}^{(d/k)-1}r^{ik}H$$ then $S$ is a subgroup of $\mathbb{Z}_{p}^{*}$ by Lemma \ref{S-T-lemma}. Thus by \ref{formula-alg-deg-circ}, $Deg(\Gamma_S)=\frac{p-1}{|S|}=\frac{p-1}{(d/k)|H|}=k<d$. Now $Orb^{(\psi_m)}(A_k)=\lbrace{r^lA_k:l=0,1,\ldots,k-1}\rbrace$ and thus $|Orb^{(\psi_m)}(A_k)|=k$, where $1\le m\le d-1$ with $gcd(m,d)>1$. Let $$T_l:=\bigcup_{i=0}^{(d/k)-1}r^{l+ik}H$$ for $0\le l\le k-1$. Then $\Gamma_S\cong\Gamma_{T_l}$ by Theorem \ref{sufficient-circ-isomorphism-theorem}. Hence $Deg(\Gamma_{T_l})=k$. Hence we can discard at least $$\sum_{\tiny\begin{array}{cc}
     k\mid d \\
     k\neq 1,d
\end{array}}k=\sigma(d)-1-d$$ many isomorphism classes of circulant graphs of order $p$ with algebraic degree less than $d$. Thus we obtain our final expression as,
$$|\mathcal{F}_{p,d}|\le\frac{1}{d}\left(\sum_{\tiny\begin{array}{cc}
     1\le m\le d-1  \\
     gcd(m,d)=1
\end{array}}{d\choose m}\right)+\frac{1}{d}(d-\varphi(d))\left(\sum_{\tiny\begin{array}{cc}
     1\le m\le d-1  \\
     gcd(m,d)>1
\end{array}}{d\choose m}\right)-(\sigma(d)-1-d)$$
$$=\frac{1}{d}(2^d-2)+\frac{1}{d}(d-\varphi(d)-1)\left(\sum_{\tiny\begin{array}{cc}
     1\le m\le d-1  \\
     gcd(m,d)>1
\end{array}}{d\choose m}\right)-(\sigma(d)-1-d).$$
It is worthy to mention that the expression obtained above will work fine as well if we consider $d$ to be prime. Thus the result holds for any $d\in\mathbb{N}$ with $d>1$ as desired. This completes the proof.\qed

If $d$ is prime, then the expression in Theorem \ref{counting-prime-ord-d-UB} will reduce down to $|\mathcal{F}_{p,d}|\le\frac{1}{d}(2^d-2)$. In the next theorem we show that this upper bound is attained when $d$ is prime. 


\begin{theorem}\label{exact-value-F_{p,d}}
Let $d$ be a prime. Then for any prime $p\in\mathcal{N}_d$, $$|\mathcal{F}_{p,d}|=\frac{1}{d}(2^d-2).$$
\end{theorem}
\pf We follow the same notations from Theorem \ref{counting-prime-ord-d-UB}. Let $d$ be a prime and let $p\in\mathcal{N}_d$ be prime. Let $\Gamma_S:=Cay(\mathbb{Z}_{p},S)$ such that $Deg(\Gamma_S)=d$. Then $S$ is a union of some $H$-orbits, i.e., $$S=\bigcup_{j\in D}r^jH$$ for some $D\subseteq\lbrace{0,1,\ldots,d-1}\rbrace$ with $1\le |D|\le d-1$, where $r$ is a primitive root in $\mathbb{Z}_{p}^{*}$. Since $d$ is prime, by Lemma \ref{group-subset-lemma-2} and \ref{formula-alg-deg-circ}, $Deg(\Gamma_S)=d$ for any choice of $D\subseteq\lbrace{0,1,\ldots,d-1}\rbrace$ with $1\le |D|\le d-1$. Thus it follows from Remark \ref{Adam_conjecture} that each of the orbits of the action $\psi_m$ actually denotes each isomorphism class of circulant graphs of order $p$ with valency $\frac{m(p-1)}{d}$ having algebraic degree exactly $d$.

Moreover, since $d$ is prime, $Fix^{(\psi_m)}(r^iH)=\emptyset$ for every $1\le i\le d-1$. By Burnside Lemma, for each $1\le m\le d-1$,
$$\#Orbits(\psi_m)=\frac{1}{d}\sum_{i=0}^{d-1}|Fix^{(\psi_m)}(r^iH)|=\frac{1}{d}\left(|Fix^{(\psi_m)}(H)|+\sum_{i=1}^{d-1}|Fix^{(\psi_m)}(r^iH)|\right)=\frac{1}{d}{d\choose m}.$$
Thus for each $1\le m\le d-1$, there are exactly $\frac{1}{d}{d\choose m}$ non-isomorphic $d$-integral circulant graphs of order $p$ with valency $\frac{m(p-1)}{d}$. Summing up the number of orbits of the actions $\psi_m$ over $m=1,\ldots,d-1$ we have the following expression which completes the proof,
$$|\mathcal{F}_{p,d}|=\#Orbits(\psi_m)_{m=1}^{d-1}=\sum_{m=1}^{d-1}\frac{1}{d}{d\choose m}=\frac{1}{d}(2^d-2).$$\qed

We end this section by illustrating Theorem \ref{exact-value-F_{p,d}} with a few examples.\\
\textbf{Example 1.} Let $d=2$ and $p=13$. Then $p\in\mathcal{N}_2$. By Theorem \ref{exact-value-F_{p,d}}, $|\mathcal{F}_{13,2}|=1$, i.e., up to isomorphism there exists a \textit{unique} $2$-integral circulant graph of order $13$. We can explicitly give the corresponding graph by [\cite{monius2022splitting},Theorem $6.1$]. Take a primitive root say, $r=2$ in $\mathbb{Z}_{13}^{*}$ and consider the subgroup $H=\langle{r^d}\rangle=\lbrace{\pm 1,\pm 3,\pm 4}\rbrace$. Then $H$ is the maximum subgroup of $\mathbb{Z}_{13}^{*}$ fixing itself. Hence the unique (up to isomorphism) $2$-integral circulant graph of order $13$ is $$\Gamma=Cay(\mathbb{Z}_{13},\lbrace{\pm 1,\pm 3,\pm 4}\rbrace).$$
\textbf{Example 2.} Let $d=3$ and $p=19$. Then $p\in\mathcal{N}_3$. By Theorem \ref{exact-value-F_{p,d}}, $|\mathcal{F}_{19,3}|=2$, i.e., up to isomorphism there exist \textit{exactly two} $3$-integral circulant graphs of order $19$. We characterize the corresponding graphs. Take a primitive root say, $r=2$ in $\mathbb{Z}_{19}^{*}$ and consider the subgroup $H=\langle{r^d}\rangle=\lbrace{\pm 1,\pm 7,\pm 8}\rbrace$. Then $\mathbb{Z}_p^{*}/H=\lbrace{H,rH,r^2H}\rbrace$. It can be seen from Theorem \ref{counting-prime-ord-d-UB} that,
$$Orb^{(\psi_1)}(H)=\lbrace{H,rH,r^2H}\rbrace,$$ $$Orb^{(\psi_2)}(H\cup rH)=\lbrace{H\cup rH,rH\cup r^2H,H\cup r^2H}\rbrace.$$
Let $S_1:=H$ and $S_2:=H\cup rH$. Then $H$ is the maximum subgroup of $\mathbb{Z}_{19}^{*}$ fixing both $S_1$ and $S_2$. Hence the non-isomorphic $3$-integral circulant graphs of order $19$ are
$$\Gamma_1=Cay(\mathbb{Z}_{19},S_1)=Cay(\mathbb{Z}_{19},\lbrace{\pm 1,\pm 7,\pm 8}\rbrace),$$
$$\Gamma_2=Cay(\mathbb{Z}_{19},S_2)=Cay(\mathbb{Z}_{19},\lbrace{\pm 1,\pm 2,\pm 3,\pm 5,\pm 7,\pm 8}\rbrace).$$
\textbf{Example 3.} Let $d=5$ and $p=11$. Then $p\in\mathcal{N}_5$. By Theorem \ref{exact-value-F_{p,d}}, $|\mathcal{F}_{11,5}|=6$, i.e., up to isomorphism there exist \textit{exactly six} $5$-integral circulant graphs of order $11$. We characterize the corresponding graphs. Take a primitive root say, $r=2$ in $\mathbb{Z}_{11}^{*}$ and consider the subgroup $H=\langle{r^d}\rangle=\lbrace{\pm 1}\rbrace$. Then $\mathbb{Z}_p^{*}/H=\lbrace{H,rH,r^2H, r^3H,r^4H}\rbrace$. It can be seen from Theorem \ref{counting-prime-ord-d-UB} that,
$$Orb^{(\psi_1)}(H)=\lbrace{H,rH,r^2H,r^3H,r^4H}\rbrace,$$ $$Orb^{(\psi_2)}(H\cup rH)=\lbrace{H\cup rH,rH\cup r^2H,r^2H\cup r^3H,r^3H\cup r^4H,H\cup r^4H}\rbrace,$$
$$Orb^{(\psi_2)}(H\cup r^2H)=\lbrace{H\cup r^2H,rH\cup r^3H,r^2H\cup r^4H,H\cup r^3H,rH\cup r^4H}\rbrace,$$
$$Orb^{(\psi_3)}(H\cup rH\cup r^2H)=\lbrace{H\cup rH\cup r^2H,rH\cup r^2H\cup r^3H,r^2H\cup r^3H\cup r^4H,}$$
$$~~~~~~~~~~~~~~~{r^3H\cup r^4H\cup H,r^4H\cup H\cup rH}\rbrace,$$
$$Orb^{(\psi_3)}(H\cup rH\cup r^3H)=\lbrace{H\cup rH\cup r^3H,rH\cup r^2H\cup r^4H,r^2H\cup r^3H\cup H,}$$
$$~~~~~~~~~~~~~~~~~~~{r^3H\cup r^4H\cup rH,r^4H\cup rH\cup r^2H}\rbrace,$$
$$Orb^{(\psi_4)}(H\cup rH\cup r^2H\cup r^3H)=\lbrace{H\cup rH\cup r^2H\cup r^3H,rH\cup r^2H\cup r^3H\cup r^4H,}$$
$$~~~~~~~~~~~~~~~~{r^2H\cup r^3H\cup r^4H\cup H,r^3H\cup r^4H\cup H\cup rH,r^4H\cup H\cup rH\cup r^2H}\rbrace.$$
Let $S_1:=H, S_2:=H\cup rH, S_3:=H\cup r^2H, S_4:=H\cup rH\cup r^2H, S_5:=H\cup rH\cup r^3H, S_6:=H\cup rH\cup r^2H\cup r^3H$. Then $H$ is the maximum subgroup of $\mathbb{Z}_{11}^{*}$ fixing $S_1,\ldots,S_6$. Hence the non-isomorphic $5$-integral circulant graphs of order $11$ are
$$\Gamma_1=Cay(\mathbb{Z}_{11},S_1)=Cay(\mathbb{Z}_{11},\lbrace{\pm 1}\rbrace),$$
$$\Gamma_2=Cay(\mathbb{Z}_{11},S_2)=Cay(\mathbb{Z}_{11},\lbrace{\pm 1,\pm 2}\rbrace),$$
$$\Gamma_3=Cay(\mathbb{Z}_{11},S_3)=Cay(\mathbb{Z}_{11},\lbrace{\pm 1,\pm 4}\rbrace),$$
$$\Gamma_4=Cay(\mathbb{Z}_{11},S_4)=Cay(\mathbb{Z}_{11},\lbrace{\pm 1,\pm 2,\pm 4}\rbrace),$$
$$\Gamma_5=Cay(\mathbb{Z}_{11},S_5)=Cay(\mathbb{Z}_{11},\lbrace{\pm 1,\pm 2,\pm 3}\rbrace),$$
$$\Gamma_6=Cay(\mathbb{Z}_{11},S_6)=Cay(\mathbb{Z}_{11},\lbrace{\pm 1,\pm 2,\pm 3,\pm 4}\rbrace).$$

\section{Conclusion and Open Issues}
As some concluding remarks we state below some of the open issues obtained so far.
\begin{itemize}
\item From Table \ref{tab:my_label}, it is observed that in most of the cases the value of $C(d)$ coincides with the upper bound $p_d$ in Corollary \ref{cor-circ-min-ord}, while they differ in some particular cases (shown in \textcolor{red}{red} font). In view of this, we pose the following question.\\
\textbf{Question:} Characterize those $d\in\mathbb{N}$ for which $C(d)=p_d$ holds.
\item In Theorem \ref{counting-d-int-circ-graph}, we hope the lower bound on the order of $\mathcal{F}_{n,d}$ can be improved by considering the connection sets suitably.
\item Finding some potential upper bound on the order of $\mathcal{F}_{n,d}$ can be considered as another topic of further work.
\item Similar to the notion $C(d)$, one can try to obtain the value of the order of the minimal $d$-integral graph (i.e., $d$-integral graph with minimum possible order). However, it seems that finding the value of this parameter in general, even providing some potential bounds bears much difficulty.
\end{itemize}

\section*{Acknowledgement}
The first author is supported by the funding of UGC [NTA Ref. No. 211610129182], Govt. of India. The second author acknowledges the funding of DST-FIST Sanction no. $SR/FST/MS-I/2019/41$ and DST-SERB-MATRICS Sanction no. $MTR/2022/000020$, Govt. of India. 

The authors are grateful to Prof. Jyrki Lahtonen, Department of Mathematics and Statistics, University of Turku, Finland for providing an important input in proving Theorem \ref{counting-prime-ord-d-alg-deg-non-iso-circ-graphs}.

\bibliographystyle{abbrv}
\bibliography{library}





\end{document}